\documentclass[12pt]{article}
\usepackage{epsfig}
\usepackage{rotating}
\title{On Pappus and Anosov Representations of
  the Modular Group}
\author{Richard Evan Schwartz \thanks{\hskip 5 pt 
Supported by 
N.S.F. Grant DMS-2505281}}

\newtheorem{theorem}{Theorem}[section]

\newtheorem{lemma}[theorem]{Lemma}

\def\startproof{{\bf {\medskip}{\noindent}Proof: }}

\def\endproof{$\spadesuit$  \newline}

\def\C{\mbox{\boldmath{$C$}}}%
\def\P{\mbox{\boldmath{$P$}}}%
\def\R{\mbox{\boldmath{$R$}}}%
\def\Z{\mbox{\boldmath{$Z$}}}%

\begin{document}

\maketitle

\begin{abstract}
  Let $X=SL_3(\R)/SO(3)$.
  Let $\cal DFR$ be the
 space of discrete 
  faithful representations of the
  modular group into
  ${\rm Isom\/}(X)$ which map the order
  $2$ generator to an isometry with a unique fixed point. In this
  paper, we prove that
  $\cal DFR$ has a component
  $\cal B$, the so-called
  Barbot component, that is homeomorphic
  to    $\R^2 \times [0,\infty)$.  The boundary of $\cal B$
  parametrizes the Pappus representations and
  the interior consists of Anosov representations.
  \end{abstract}

      \section{Introduction}

Let  $X=SL_3(\R)/SO(3)$.  This is a prototypical
higher rank symmetric space.
In this paper we completely characterize one
connected component of the moduli space
$\mathcal DFR$ of conjugacy classes of
discrete and faithful
representations of the modular group
$\Z/3* \Z/2$ into
${\rm Isom\/}(X)$ which map
the order $2$ elements to isometries
having a unique fixed point in $X$.

The {\it Pappus representations\/}
  are a $2$-parameter subfamily
  of $\mathcal DFR$ which I constructed in
  my 1993 paper
  [{\bf S0\/}] and then revisited
  in my recent paper [{\bf S1\/}].
  These groups exhibit many features
  that, much later and more generally,
  appeared in
   higher Teichmuller Theory, e.g. in [{\bf Lab\/}],
  [{\bf GW\/}], [{\bf Bar\/}],
  [{\bf BCLS\/}], and [{\bf KL\/}].
  
The Pappus representations
are
  nowadays classified as {\it relatively Anosov groups in the
    Barbot component\/}.
  This point of view is exposited in
  [{\bf BLV\/}] and [{\bf KL\/}].
  Let ${\mathcal P\/} \subset \mathcal DFR$ denote
  the subset consisting of Pappus modular
  group representations.  The connected component
  $\mathcal B$ of $\mathcal DFR$ containing
  $\mathcal P$ is called the {\it Barbot component\/}.
  It is partially understood
  thanks to  [{\bf S0\/}] and [{\bf BLV\/}].
  
  In [{\bf BLV\/}], T. Barbot,
  G.-S. Lee, and V. P. Valerio build
    on [{\bf S0\/}] and construct
  a $3$-parameter family of Anosov
  representations which are defined
  in terms of modified operations on
  marked boxes.  Using their
  {\it morphed marked boxes\/} (my terminology)
  they construct a $4$-parameter family of
  representations of $\Z/3 * \Z/3$ into
  $SL_3(\R)$, all of which are Anosov.
  They show that a subset $\mathcal A$ of these
  extend to Anosov representations of
  the modular group $\Z/3 * \Z/2$.
  They use an implicit
  function argument to show the existence of
  $\mathcal A$.
  Here I will
  complete the analysis.
    
\begin{theorem}
  \label{one}
  $\cal DFR$ has a connected component $\cal B$ which is homeomorphic
    to $\R^2 \times [0,\infty)$.   The representations
  corresponding to $\R^2 \times \{0\}$ are the
  Pappus modular groups, and the representations
  corresponding to $\R^2 \times (0,\infty)$ are Anosov
  representations.
\end{theorem}

Here is more information about $\cal B$.
We first define a larger representation space
$\cal R$ of all sufficiently generic
representations of the group
$\Z/3 * \Z/2$ into ${\rm Isom\/}(X)$ which
map the order $2$ element to an isometry
having a unique fixed point in $X$.
Here, {\it sufficiently generic\/} means
that the fixed point of the order $2$ element
is not contained in the fixed set of the order $3$
element.
We show that $\cal R$ is homeomorphic
to $\R^3-\{(0,0,0)\}$ and that
$\cal P$ is a properly embedded plane
in $\cal R$, smooth away from one point.
We show that $\cal B$ is
the closure of the component
of ${\cal R\/}-{\cal P\/}$ that does not contain the origin.
The boundary
of $\cal B$ is exactly $\cal P$.

Our proof of these results
involves extending the analysis in [{\bf BLV\/}] to
fully work out the set $\mathcal A$.
Our trick is to replace the transcendental
parametrization in [{\bf BLV\/}] with a rational parametrization
and then subject the resulting formulas to computer algebra.

This paper is a short version of a longer paper [{\bf S2\/}] I wrote,
which goes much more deeply into the structure of groups in
$\cal B$.  I thought that Theorem \ref{one} would
be interesting on its own, and that
a shorter paper just dealing with this one result would be
more accessible.  Also, one step in [{\bf S2\/}] has an error.
In brief, I analyzed the duality curve associated to parameters
$b \in (0,1)$ when I meant to analyze the duality curve associated
to parameters $b \in (1,\infty)$.  The analysis is essentially the
same
but certainly this needs to be fixed.  So, with this paper, I take
the opportunity to actually analyze the correct parameter set.

This paper is organized as follows.
    \begin{itemize}
\item  In \S 2, I give some background material
  about the space $X$.
\item In \S 3 I give a
  topological analysis of the representation
  space $\mathcal R$.

\item In \S 4, I give an exposition of the Pappus modular
  groups that is similar to what I gave in [{\bf S1\/}].
  
\item In \S 5,  I recall the work in [{\bf BLV\/}] and then
  recast their construction in more algebraic terms.

  \item In \S 6, I use
  algebraic methods to completely
  analyze the set $\mathcal A$ of Anosov representations
  studied in [{\bf BLV\/}].

\end{itemize}

The calculations in this paper are done in Mathematica.
One can download the Mathematica code for this paper
from \newline
{\bf http://www.math.brown.edu/$\sim$res/PappusCalcs.TAR\/}
\newline
If you load the files (one at a time) into Mathematica, they will
print out essentially every relevant formula;  you can check
these calculations against the text of the paper.
I think that this paper is best read with these Mathematica files
in hand.

I  thank
Martin Bridgeman,
Bill Goldman,  Tom Goodwillie,
Sean Lawton,
Joaquin Lejtreger,
Joaquin Lema, Dan Margalit,  Max
Riestenberg, Dennis
Sullivan, and Anna Wienhard for
interesting and helpful conversations.
I especially thank Martin and the
two Joaquins for recently rekindling
my interest in this subject.

      \section{Geometric Preliminaries}

\subsection{Projective Geometry}

The projective plane $\P$ is the set of
$1$-dimensional subspaces of $\R^3$. 
The dual plane $\P^*$ is the set of
$2$-dimensional subspaces of $\R^3$. 
We represent points in $\P$ in homogeneous
coordinates.
When $c \not =0$, the
element $[a:b:c]$ represents
$(a/c,b/c)$ in the {\it affine patch\/}.
The affine patch is a copy of
$\R^2$ sitting inside $\P$.
Any inner product $g$ on $\R^3$ determines
an analytic diffeomorphism from
$\P$ to $\P^*$.  The point of $\P$ represented by
a line $L$ through the origin is mapped to
the point of $\P^*$ represented by the
$g$-perpendicular complement. We call
such an element an {\it elliptic polarity\/}. When
$g$ is the usual dot product, we call this
map the {\it standard elliptic polarity\/} and
denote it by $\Delta$.
A {\it duality\/} is the
composition of a projective transformation
with the standard elliptic polarity.

We represent projective transformations by
elements of $GL_3(\R)$, the group of
invertible $3 \times 3$ matrices.
Every such element has a unique representation
in $SL_3(\R)$, but in the interest of getting
rational expressions we will sometimes
refrain from scaling.  We can represent a duality
as $\delta \circ M$ where $M \in PSL_3(\R)$.

Given any matrix
$m \in GL_3(\R)$, the quantity
\begin{equation}
  \label{TAU}
  \tau(m)=\frac{{\rm tr\/}^3(m)}{{\rm det\/}(m)}.
\end{equation}
is independent of the scaling of $m$ and also is
a $GL_3(\R)$-conjugacy invariant.

\subsection{The Symmetric Space}
\label{adjoint}

My article [{\bf S1\/}] has detailed information about the symmetric
space
\begin{equation}
  X=SL_3(\R)/SO(3).
\end{equation}
Here I give an abbreviated account.
$X$ is the space of unit determinant positive definite symmetric
$3 \times 3$ matrices.  The linear group  $SL_3(\R)$ acts
isometrically on $X$.   The action is given by
\begin{equation}
  T(M)=T^* \circ M \circ T^{-1}
\end{equation}
Here $T^*$ is the inverse-transpose of $T$.
This action has the following interpretation.
If $E$ is the unit ball for $M$ then
$T(E)$ is the unit ball for $T(M)$.

$X$ has a natural origin, namely the point $O$ given by
the identity matrix.  We let
$M(a,b,c)$ denote the diagonal matrix with
entries $a,b,c$.  Here we have $abc=1$.
Thus, $O=M(1,1,1)$.
The space $X$ has a canonical Riemannian metric
with respect to which $SL_3(\R)$ acts isometrically.
The distance between $O$ and $M(a,b,c)$ is given by
\begin{equation}
  \sqrt{\log^2(a)+\log^2(b)+\log^2(c)}.
\end{equation}
The rest of the metric can be deduced from
symmetry.

The standard elliptic polarity $\Delta$ induces an action on $X$, and
the action is given by:
\begin{equation}
  \Delta(M)=M^{-1}.
\end{equation}
This map is an involution which reverses all the geodesics
through $O$.  The point $O$ is the only fixed point of $\Delta$.
The isometry 
group ${\rm Isom\/}(X)$ is generated by
$\Delta$ and $SL_3(\R)$.

\subsection{The Adjoint Action}

Let $T_O(X)$ denote the tangent space to $X$ at $O$.
This is the space of trace $0$ symmetric matrices.
The subgroup $SO(3)$ acts on $T_O(X)$ by the
{\it adjoint representation\/}:
\begin{equation}
  g: M \to gMg^{-1}.
\end{equation}
The reader might worry that we should really
use $(g^{-1})^t$ in place of $g$ but fortunately
$g=(g^{-1})^t$ when $g \in SO(3)$.
Also, we mention that technically we are talking
about the {\it restriction\/} of the adjoint action
to a maximal compact subgroup of $SL_3(\R)$.

For purposes that will be made clear in the
next section we wish to consider the adjoint
action of the matrices
\begin{equation}
  \label{rotation}
  \left[\matrix{\cos(\theta)&\sin(\theta)&0 \cr
      -\sin(\theta)&\cos(\theta)&0 \cr 0&0&1}\right]: \hskip 30 pt
  \left[\matrix{a&b&c \cr b& -a & d \cr c& d & 0}\right] \to
  \left[\matrix{a'&b'&c' \cr b'& -a' & d'\cr c'& d' & 0}\right].
\end{equation}
We calculate that
  $$a' =  a \cos(2 \theta) + b \sin(2\theta),
  \hskip 30 pt b'=-a\sin(2\theta)+b \cos(2\theta),$$
  $$
    c'=c\cos(\theta) + d\sin(\theta), \hskip 30 pt
   d'= -c \sin(\theta) + d\cos(\theta).
    $$
    This action looks nicer if we identify the matrix in
    Equation \ref{rotation} with the unit complex number
    $u=\exp(i\theta)$ and
    $\R^4$ with $\C^2$ under the identification
    \begin{equation}
      (a,b) \to z=a+bi, \hskip 30 pt (c,d) \to w=c+di.
    \end{equation}
    The action is then given by
    \begin{equation}
      \label{act1}
      u: (z,w) \to (u^2 z,uw).
    \end{equation}

    Here is the geometric significance of the matrices
    on the right side of Equation \ref{rotation}.
    They are all orthogonal to the tangent vector given
    by the matrix ${\rm diag\/}(-1,-1,2)$.   This matrix
    (considered as a tangent vector)
    is in turn tangent to the geodesic in $X$
    through the origin that limits at one end on the
    point of $\P$ named by the origin and at the other
    end on the point of $\P^*$ named by the line at
    infinity.   We let $V_0 \cong \C^2$ be the
    vector space of such matrices.

    Our action acts in a rather special way on the subspaces
    $\C \times \{0\}$ and $\{0\} \times \C$.   The former
    subspace is the tangent space to matrices which have
    block form with a $2 \times 2$ matrix in the upper
    left corner and a nonzero entry in the lower right corner.
    The latter subspace corresponds to matrices which
    stabilize the unit circle in the affine patch.  These
    two subspaces will correspond to representations which,
    respectively, preserve a projective line and a conic section.
    The former arise for us and the latter do not.

      \section{The Representation Space}
\label{bigrep}

See e.g. [{\bf L\/}] and [{\bf FL\/}] for
topics related to the material here.

The modular group $G=\Z/2*\Z/3$ is generated by
$\sigma_2$ and $\sigma_3$, elements of order $2$ and $3$.
We consider representations of $G$ into ${\rm Isom\/}(X)$ such that
\begin{itemize}
\item $\rho(\sigma_2)$ is an elliptic polarity.
\item  $\rho(\sigma_3)$ is
  the matrix in Equation \ref{rotation} for
$\theta=2\pi/3$.   This matrix acts on
$\P$ as an order $3$ rotation of the affine patch.
\item The fixed point of $\rho(\sigma_2)$ does not
  lie in the fixed point set of $\rho(\sigma_3)$.
\end{itemize}
The fixed point set
of $\rho(\sigma_3)$ is the geodesic $\gamma$
consisting of standard ellipsoids $E(a,a,a^{-2})$.
We call these representations
{\it normalized\/}.
We consider two representations to be
the same if
they are conjugate in ${\rm Isom\/}(X)$.

Let $\cal R$ denote the space of all normalized
representations, modulo conjugacy.
We define the distance between two elements
$[\rho_1], [\rho_2] \in \cal R$
to be the minimal $D$ such that there are
two normalized representatives
$\rho_1$ and $\rho_2$ such that
the fixed point sets of $\rho_1(\sigma_2)$
and $\rho_2(\sigma_2)$ are $D$ apart in $X$.

\begin{theorem}
  \label{trace}
  $\cal R$ is homeomorphic to $\R^3-\{(0,0,0)\}$
  and is a smooth manifold away
    from the two curves, one corresponding to
    line-preserving representations and one corresponding
    to conic-preserving representations.
    The trace of  $\rho(\sigma_3\sigma_2\sigma_3\sigma_2)$
    is a smooth function away from the two special curves.
  \end{theorem}

  \noindent
  {\bf Remark:\/}
  The statement about $\sigma_3\sigma_2\sigma_3\sigma_2$
  holds more generally for any word in the group, but we only
  care about this one word.
  \newline

The rest of the chapter is devoted to proving this result.
Recall that $\gamma$ is the (singular) geodesic
fixed by $\rho(\sigma_3)$ for all $\rho$.
For each $p \in \gamma$ we let $V_p$ be
the subspace of the tangent space $T_p(X)$ which
is orthogonal to $\gamma$.   Let $X_p$ denote the
image of $V_p$ under the exponential map.
We call $X_p$ an {\it orthogonal cut\/}.
The orthogonal cuts are diffeomorphic to
$\R^4$.

\begin{lemma}
  The space $X$ is foliated by the orthogonal cuts.
  \end{lemma}

  \startproof
  Every point $q \in X-\gamma$ lies in the orthogonal
  cut containing the geodesic connecting $q$ to the
  point on $\gamma$ nearest $q$.  Given this fact, we
  just have to show that two orthogonal cuts are
  disjoint. If not, we can find a geodesic triangle in
  $X$ with $2$ right angles.   But this is impossible
  in a space like $X$, which has non-positive
  sectional curvature.
  \endproof

  Let $S_{\gamma}$ denote the stabilizer of $\gamma$
  in ${\rm Isom\/}(X)$.
      Using the action of $S_{\gamma}$ we can normalize so that the
    fixed point of $\rho(\sigma_2)$ lies in the orthogonal cut $X_0$
    through the origin.  The reason this is possible is that
    $S_{\gamma}$ acts transitively on $\gamma$ and hence
    acts transitively on the set of orthogonal cuts.
    Let $S_{\gamma}^0 \subset S_{\gamma}$ be the subgroup which stabilizes
    $X_0$.  This subgroup is generated by rotations,
    as in Equation \ref{rotation}, and the standard
    polarity.  The rotations act on $V_0 \cong \C^2$ as
    in Equation \ref{act1}, and the polarity acts as
    $\Delta(z,w)=(-z,-w)$.

    Using the inverse exponential map, a diffeomorphism,
    we identify $X_0$ with
    the $4$-dimensional subspace $V_0 \cong \C^2$ discussed in the
    previous section.   So, the quotient we want is
    \begin{equation}
      \label{QU}
      (\C^2-(0,0))/S_{\gamma}^0.
    \end{equation}
    The action of $S_{\gamma}^0$ preserves the standard polar
    coordinate system in $\C^2 \cong \R^4$, so the quotient
    we seek is just the cone (minus the origin) over
    $S^3/S_{\gamma}^0$.  We now have a standard topological
    problem.

    The quotient $S^3/S_{\gamma}^0$ is homeomorphic to
    $S^2$, and has a smooth structure away from the points
    corresponding to the circles $\{z=0\}$ and $\{w=0\}$.
    Here we recall the construction.
    Let $S^3_*$ denote the space obtained by removing these
    two circles.  This space is foliated by Clifford tori satisfying
    the equations
    $|z|/|w| = {\rm const\/}$, and the action of $S_{\gamma}^0$ preserves this
    foliation.  The quotient $S^3_*/S_{\gamma}^0$ is diffeomorphic to
    the product $(T/S_{\gamma}^0) \times (0,\infty)$ where $T$ is
    the central Clifford torus $|z|=|w|$.   The quotient
    $T/S_{\gamma}^0$ is diffeomorphic to a circle.  Hence
    $S^3_*/S_{\gamma}$ is homeomorphic to a cylinder.
    But then $S^3/S_{\gamma}^0$ is the two-point compactification
    of this smooth cylinder, a topological sphere.

    Taking the cone, we see that the quotient in
    Equation \ref{QU} is a smooth manifold away
    from the curves coming from the cones over the
    two special points of $S^3/S_{\gamma}^0$.
    This gives us everything in Theorem \ref{trace}
    except the statement about the trace.
    
    The trace of $\rho(\sigma_3\sigma_2\sigma_3\sigma_2)$ 
    is a polynomial function on the matrix entries of
    $\rho(\sigma_2)$ and $\rho(\sigma_3)$.
    (Here we represent the polarity $\rho(\sigma_2)$ as a matrix $M$
    such that $\rho(\sigma_3)=\Delta \circ M$.)
    When we construct a local coordinate chart for the smooth
    subset of the quotient in Equation \ref{QU} what we do
    is take a small and smooth  cross section to the circle foliation given
    by the action in Equation \ref{act1}.  The trace of our given word
    restricts to a smooth function on this cross section.   Hence
    ${\rm tr\/}(\rho(\sigma_3\sigma_2\sigma_3\sigma_2))$ is a smooth
function on    the smooth part of $\cal R$.

\section{The Pappus Modular Groups}

\subsection{Basic Definitions}

In this chapter we recall the Pappus
modular group representations.
Our exposition follows [{\bf S1\/}], though ultimately
the material goes back to  [{\bf S0\/}].  The paper
[{\bf BLV\/}] also has an exposition that is like [{\bf S0\/}].
\newline
\newline
{\bf Convex Marked Boxes:\/}
A {\it convex marked box\/} is a convex quadrilateral
in $\P$ together with a distinguished point in the
interior of one side and a distinguished point in the
interior of an opposite side.  We call one of the
points the {\it top\/} point and the other one the
{\it bottom\/} point.  Correspondingly we call the
edges containing these points the {\it top edge\/}
and the {\it bottom edge\/}.  Finally, we say that
the {\it top flag\/} is the flag $(p,\ell)$ where $p$
is the top point and $\ell$ is the line extending the
top edge.  We define the {\it bottom flag\/} similarly.
\newline
\newline
{\bf Operations on Marked Boxes:\/}
There are $3$ operations we can perform on marked
boxes, and we call them $t,b,i$.  Figure 4.1
shows how they act.

 \begin{center}
\resizebox{!}{2in}{\includegraphics{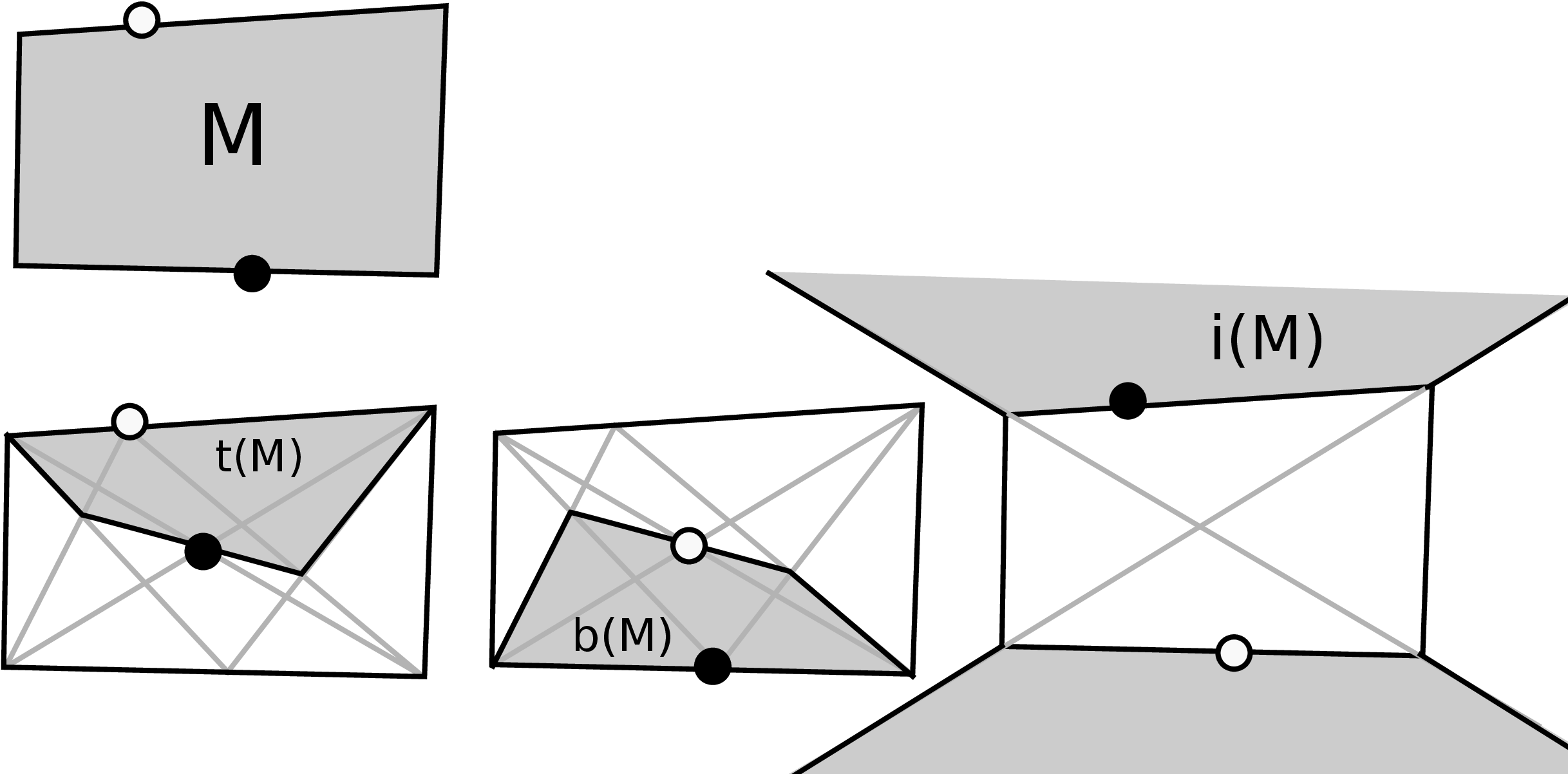}}
\newline
{\bf Figure 4.1:\/} The three operations on marked boxes
\end{center}

These operations satisfy the relations
\begin{equation}
  i^2=I. \hskip 20 pt
  tit=b, \hskip 20 pt bib=t, \hskip 20 pt
  tibi=I, \hskip 20 pt biti=I.
\end{equation}
Here $I$ is the identity.  As a consequence of these relations,
and the nesting of the marked boxes, the group of operations
isomorphic to the modular group.    The explicit generators are
(say) $i$ and $ti$.  
We let $\cal M$ be the orbit of a marked box
under the action of this group.
\newline
\newline
{\bf Order Three Symmetries of the Orbit:\/}
Given a marked box $M \in \cal M$ there is an order $3$ projective
transformation $T_M$ which has the orbit
$$i(M) \to t(M) \to b(M).$$
This accounts for the order
$3$ elements of the Pappus modular groups.
\newline
\newline
{\bf Order Two  Symmetries of the Orbit:\/}
There is also an elliptic polarity which, in a certain sense,
swaps $M$ and $i(M)$.  To make sense of this, we have
to recall the notion of a {\it doppelganger\/} defined in
[{\bf S1\/}].

 \begin{center}
\resizebox{!}{1.62in}{\includegraphics{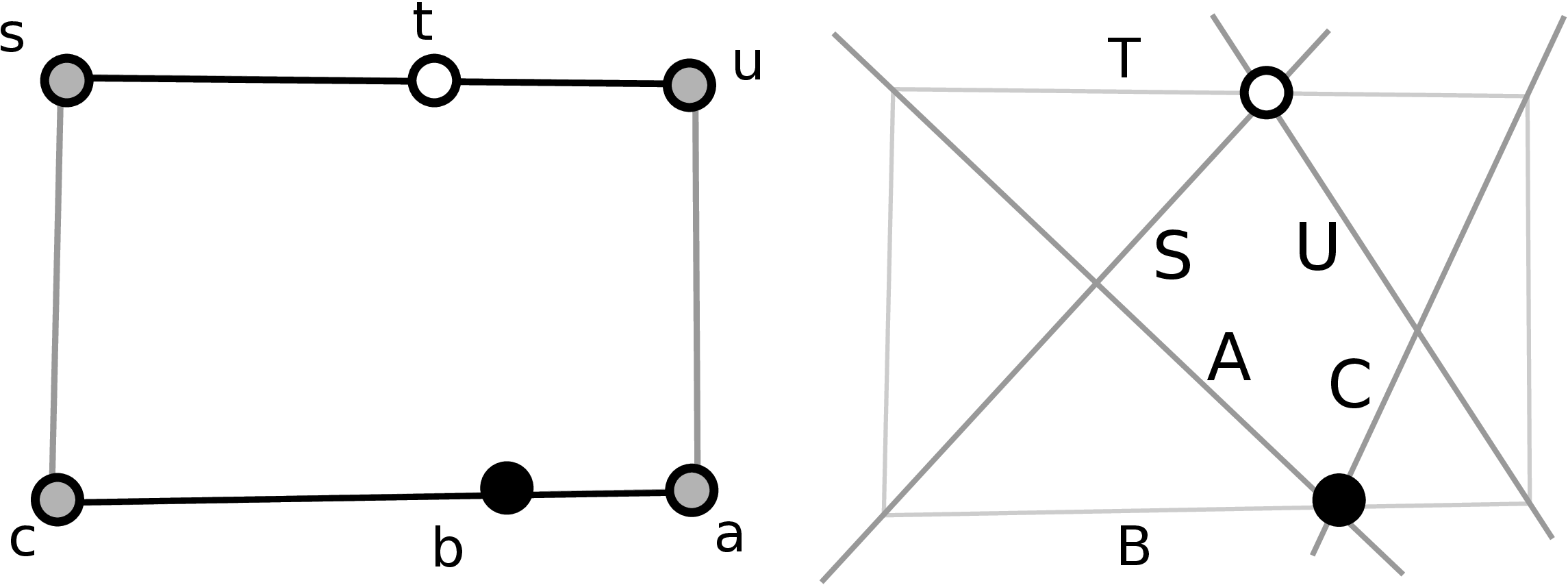}}
\newline
{\bf Figure 4.2\/} A convex marked box and its doppelganger
\end{center}

The $6$-tuple $(s,t,u,a,b,c)$ shown on the left side of Figure 4.2
encodes the marked box $M$.  Here $t$ and $b$ are respectively
the top and bottom points of $M$.  The corresponding $6$-tuple of lines
$(S,T,U,A,B,C)$, which is defined entirely in terms of $M$,
encodes a convex marked box $M^*$ in $\P^*$.    We can repeat
the operation and we get $M^{**}=M$.    It turns out that the
$i,b,t$ operations commute with the doppelganger operation and
we can think of our orbit $\cal M$ as an orbit of pairs of the
form $(M,M^*)$.  We call such a pair an {\it enhanced convex marked
  box\/}.

We showed in [{\bf S1\/}] that there is an elliptic polarity
$\delta_M$ 
that swaps $M$ and $(i(M))^*$, and simultaneously
swaps $M^*$ and $i(M)$.

We showed in [{\bf S0\/}] that the Pappus groups define
discrete and faithful representations of the modular group.
Each one is a point of the representation space
$\cal R$.

\subsection{Formulas}
\label{FORMULA}

Referring to Figure 4.2, we can normalize a marked box by
a projective transformation so that its vertices are
as in Figure 4.3.
   \begin{center}
     \resizebox{!}{1.4in}{\includegraphics{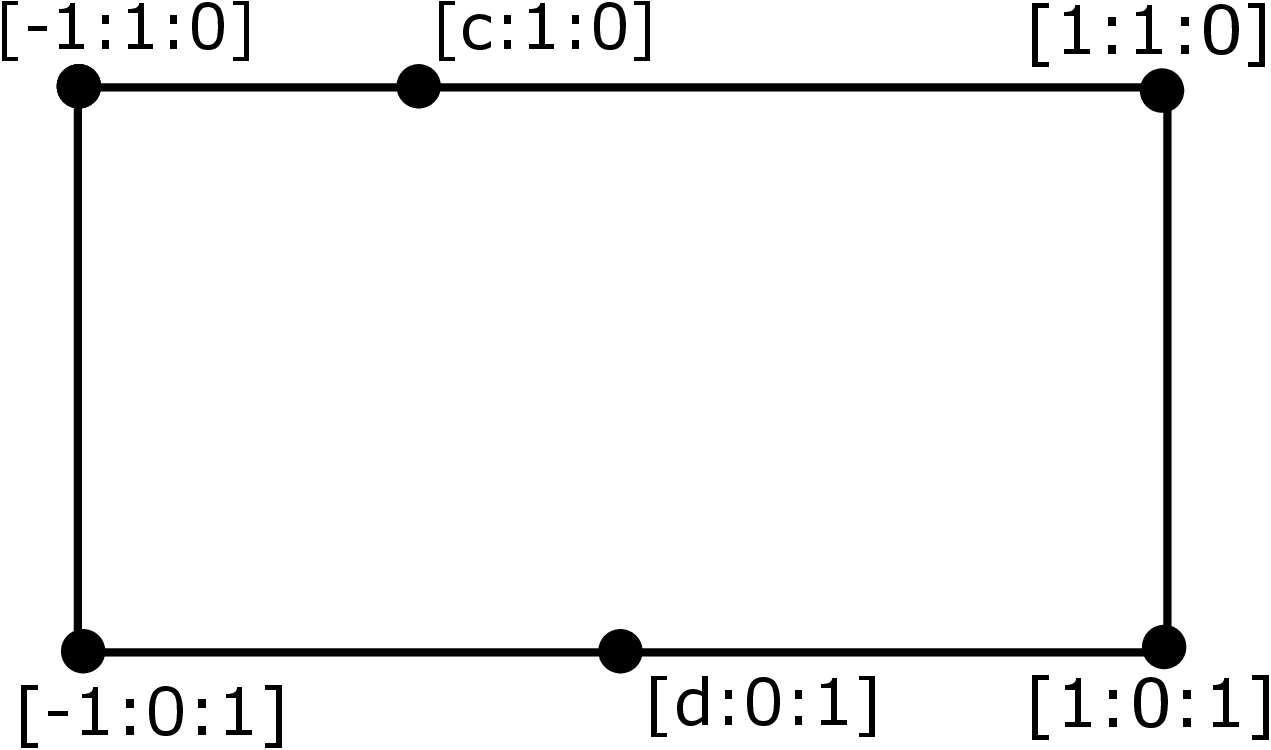}}
    \newline
    {\bf Figure 4.3\/}  The initial box
  \end{center}
  This normalization is different than the
  one in [{\bf S1\/}], but it matches [{\bf BLV\/}].

  For our purposes, we just need to know the
  index $2$ subgroup of the Pappus group that
  lies in $PGL_3(\R)$. This subgroup is generated
  by the following matrices:

$$
r_1=
\frac{1}{(1-c^2)(1-d^2)}
\left(
\begin{array}{ccc}
cd-1 & c(1-cd) & d-c \\[4pt]
d-c & 1-cd & cd-1 \\[4pt]
0 & 1-c^2 & 0
\end{array}
\right).
$$

$$r_2=
\left(
\begin{array}{ccc}
-1-cd & c+d & d(1+cd) \\
0 & 0 & d^2-1 \\
-c-d & 1+cd & 1+cd
\end{array}
\right).
$$

These matrices do not have unit determinant, but their
product does.  The product $r_1r_2$ has trace $-1$
and is parabolic. Both $r_1^3$ and $r_2^3$ are a constant
times the diagonal matrix.
Acting projectively, $r_1$ and $r_2$ respectively have the action
\begin{equation}
  \label{action}
  i(M_{c,d}) \to t(M_{c,d}) \to b(M_{c,d}), \hskip 30 pt
  M_{c,d} \to ti(M_{c,d}) \to bi(M_{c,d}).
\end{equation}
In terms of representations, we have
$r_1=\rho(\sigma_3)$
and $r_2=\rho(\sigma_2\sigma_3\sigma_2)$.
So, $r_1r_2$ is the word considered in
Theorem \ref{trace}.

Here are some additional formulas

\begin{equation}
  \label{limit}
  \tau(r_1r_2^2) =
\frac{64}{(1-c^2)(1-d^2)^2}.
\hskip 30 pt
  \tau(r_1^2r_2) =
\frac{64}{(1-c^2)^2(1-d^2)}.
\end{equation}

\begin{equation}
  \label{comm}
  {\rm tr\/}[r_2,r_1] - {\rm
    tr\/}[r_1,r_2]=\frac{16cd}{(1-c^2)(1-d^2)}.
\end{equation}
Here $\tau$ is as in Equation \ref{TAU}
and $[r_1,r_2]=r_1r_2r_1^2r_2^2$ is the commutator
of $r_1$ and $r_2$.

\subsection{The Space of Pappus Representations}
\label{papspace}

Let $\theta_4$ denote the order
$4$ rotation of $(-1,1)^2$ about $(0,0)$.

\begin{lemma}
Two pappus representations are
conjugate in ${\rm Isom\/}(X)$ if and only if
they are in the same $\theta_4$-orbit.
\end{lemma}

\startproof
We gave a geometric proof in [{\bf S1\/}].  Here
we give an algebraic proof.
Looking at Equations \ref{limit} and
\ref{comm} we see that the the $SL_3(\R)$ conjugacy
class of the representation determines
$(c,d)$ up to sign.  Hence
Hence the points $(c_1,d_1)$ and
$(c_2,d_2)$ determine
$SL_3(\R)$-conjugate representations
only if $(c_1,d_1)=\pm (c_2,d_2).$

The action of dualities on traces is trickier to understand
but it suffices to see what happens for the duality of our
choice. All other ones have the same action.  We choose to
look at the duality $\rho(\sigma_2)$.
The duality $\rho(\sigma_2)$ conjugates $(r_1,r_2)$ to $(r_2,r_1)$
and the roles of $c$ and $d$ switch.  So, if
$(c_1,d_1)$ and $(c_2,d_2)$ give representations which
are conjugate by a duality we have
$c_1^2=d_2^2$ and $d_1^2=c_1^2$.
When we swap $r_1$ and $r_2$, the sign
in Equation \ref{comm} changes, so we have
$c_1d_1=-c_2d_1$.
Putting everything together, we see that
$(c_1,d_1)$ and $(c_2,d_2)$ give representations
that are conjugate in ${\rm Isom\/}(X)$ only if
they lie in the same $\theta_4$-orbit.

For the converse, we note that
$(c_1,d_1)$ and one of the two choices
$\pm  (d_1,-c_2)$ give conjugate
representations because they are conjugate by
$\rho(\sigma_2)$.   Finally, the geometric operation
of reflecting our initial box in the $Y$-axis
conjugates the representation given
by $(c,d)$ to the one given by $(-c,-d)$.  So,
points in the same $\theta_4$-orbit give representations
that are conjugate in ${\rm Isom\/}(X)$.
\endproof

Hence,
the space of Pappus representations modulo
conjugacy is given by the quotient
\begin{equation}
  \label{quotient}
  {\cal C\/}=(-1,1)^2/\theta_4.
\end{equation}
Here ${\cal C\/}$ is an open cone that
is homeomorphic to $\R^2$.

We have a map $f: {\cal C\/} \to {\cal R\/}$, which
maps the Pappus representation parameterized
by $[(c,d)] \in \cal C$ to the point in $\cal R$ that names it.

\begin{lemma}
  \label{proper}
  $f$ is a proper map.
\end{lemma}

\startproof
Evidently the traces in
Equation \ref{limit} tend to
$\infty$
If either $|c| \to 1$ or $|d| \to 1$.
This means that the corresponding representation
exits every compact subset of $\cal R$ if the
corresponding parameter $[(c,d)]$ exits every
compact subset of $\cal C$.
\endproof

Because $f$ is a proper map, smooth away from
the totally symmetric point of $\cal C$, the image
${\cal P\/}=f({\cal C\/})$
separates ${\cal R\/}$ into two components.
This is a consequence of the Jordan
Separation Theorem.
It turns out that away from the totally
symmetric point $\cal P$ is a smooth
embedded $2$-manifold.   The idea is that
locally $\cal P$ is the level set of the
function ${\rm tr\/}(r_1r_2)=-1$.
See the end of \S \ref{DUALITY} for more discussion.

      \section{The Anosov Picture}

\subsection{Morphing Marked Boxes}

   The construction in [{\bf BLV\/}] builds off the marked box
   construction from [{\bf S0\/}].   Here we recall the
   constructions in [{\bf BLV\/}].
   Barbot, Lee, and Valerio identify a certain operation
   $\sigma_{\delta,\epsilon}$ which modifies a marked box
   by a projective transformation.  Here $\delta$ and
   $\epsilon$ are real parameters.  This is really
   an operation on convex quadrilaterals; the distinguished
   top and bottom points just go along for the ride.
   Figure 5.1 shows the image of the unit square under
   $\sigma_{-1/5,-1/5}$.

    They define their operation in a way that forces it to be
    projectively
    natural.  Given a marked box $M$ they let $T_M$ be a projective
    transformation so that $T_M(M)$ has vertices
    \begin{equation}
      \label{vert}
    [-1:1:0], \hskip 30 pt [1:1:0], \hskip 30 pt
    [1:0:1], \hskip 30 pt [-1:0:1].
  \end{equation}
  We call this particular quadrilateral $M_0$.
    These points are listed so that they go cyclically around the
    boundary of the convex quad.   The first two vertices are on the
    top edge and the last two vertices are on the bottom edge.
    $T_M$ is unique up to an order $2$ symmetry.
    Next, they introduce the projective transformation given by
    \begin{equation}
      \label{SIGMA}
      \Sigma_{\delta,\epsilon}=\left[\matrix{1&0&0 \cr 0  & e^{-\delta}
          \cosh(\epsilon) & -\sinh(\epsilon) \cr
          0&-\sinh(\epsilon) & e^{\delta} \cosh(\epsilon)}\right].
    \end{equation}
    Finally, they define
    \begin{equation}
      \sigma(M)=T_M^{-1} \circ \Sigma \circ T_M.
    \end{equation}
    See [{\bf BLV\/}, \S 7.1].  Let's call this  {\it box morphing\/}.

   \begin{center}
     \resizebox{!}{2in}{\includegraphics{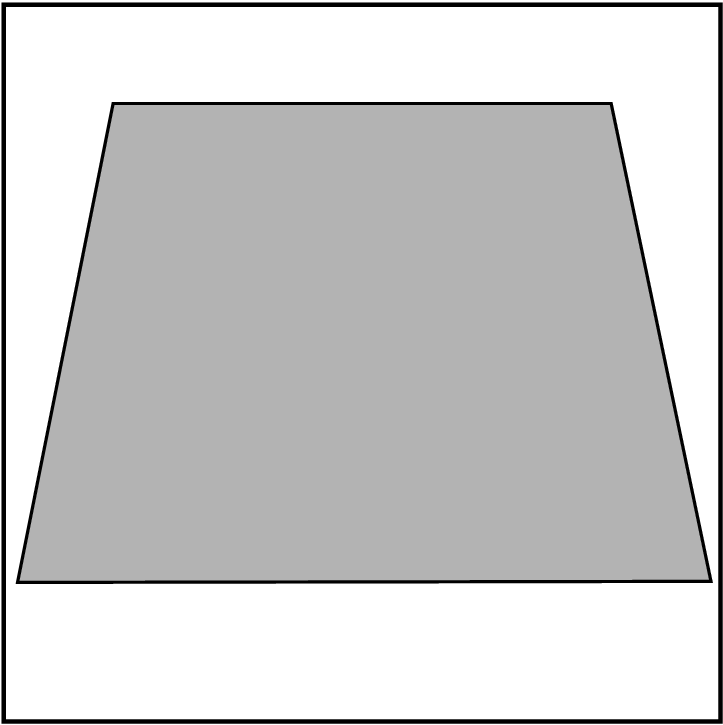}}
    \newline
    {\bf Figure 5.1\/}  The unit square morphed by
    $\sigma_{-1/5,-1/5}$.
    \end{center}

    Now we come to our main idea.
     We replace the transcendental functions in
    Equation \ref{SIGMA} with rational functions.  We set
    \begin{equation}
      a=e^{\delta}, \hskip 30 pt
      \sinh(\epsilon)=\frac{1-b^2}{2b}, \hskip 30 pt
      \cosh(\epsilon)=\frac{1+b^2}{2b}.
    \end{equation}
    Here $(a,b) \in (0,\infty)^2$.
    These are rational parametrizations of these transcendental
    functions.  We now define

    \begin{equation}
      \Sigma_{a,b}=\left[\matrix{1&0&0 \cr 0 & \frac{(1+b^2)}{2ab} &
          \frac{-1+b^2}{2b} \cr
          0 & \frac{-1+b^2}{2b} & \frac{a(1+b^2)}{2b}} \right].
    \end{equation}

We define the set $\Theta$ of {\it good parameters\/}
$(a,b)$ to be those parameters for which
$\Sigma_{a,b}(M_0)$ is contained in the interior
of  $M_0$.

\begin{lemma}
  \label{GOOD}
  $\Theta$ has the following description.
  When $b \in (1,1+\sqrt 2)$ we have
  $$\frac{1+2b-b^2}{b^2+1}<a< \frac{b^2+1}{1+2b-b^2}.$$
  When $b \in [1+\sqrt 2,\infty)$ we have $a \in (0,\infty)$.
  \end{lemma}

  \startproof
For $(a,b) \in \Theta$ it is certainly
necessary that $\Sigma_{a,b}$ maps
each of the vertices of $M_0$ into
the interior of $M_0$.  However,
this is not quite sufficient.  We also need
to check that $\Sigma_{a,b}$ maps one point
of each edge of $M_0$ into the interior
of $M_0$.
  The constraints just stated define a connected
  subset of the space $(0,\infty) \times (1,\infty)$, and the image
  of $\Sigma_{a,b}(M_0)$ varies continuously with
  the parameters $(a,b)$.  For this reason, it suffices to
  prove that the conditions of the lemma describe when
  $\Sigma_{a,b}$ maps the vertices of $M_0$ into
  $M_0$.

We will use homogeneous coordinates for our calculations.
We first note a symmetry: $\Sigma_{a,b}$ commutes with
reflection in the $y$-axis and $M_0$ (as a convex
quadrilateral) is symmetric with respect to reflection
in the $y$-axis.  For this reason, it suffices to check
the two right vertices of $M_0$, namely
$[1:1:0]$ and $[1:0:1]$.   The images of these points
under $\Sigma_{a,b}$, in the affine patch, are respectively
\begin{equation}
  \bigg(\frac{2b}{b^2-1},\frac{b^2+1}{a(b^2-1)}\bigg),
  \hskip 30pt
  \bigg(\frac{2b}{a(b^2+1)},\frac{b^2-1}{a(b^2+1)}\bigg).
\end{equation}
These points lie in the positive quadrant.  To
lie in $M_0$ they must lie above the
line $y=x-1$.  That is, they must satisfy $y>x-1$.

Applying this to our two points, we get the constraints
$$(b^2+1)+a(b^2-2b-1)>0, \hskip 30 pt
a(b^2+1)+(b^2-2b-1)>0.$$
When $b \geq 1+\sqrt 2$ all terms are positive, and
any $a>0$ works.
When $b \in (1,1+\sqrt 2)$, we have the constraints
advertised in the lemma.
\endproof

Figure 5.2 shows part of the region $\Theta$.
One should compare Figure 11 in [{\bf BLV\/}].
The region below the line $b=1$ is not relevant
to Theorem \ref{one}.
For $b>1$ the two constraint curves
open up monotonically on either side.

   \begin{center}
     \resizebox{!}{3.2in}{\includegraphics{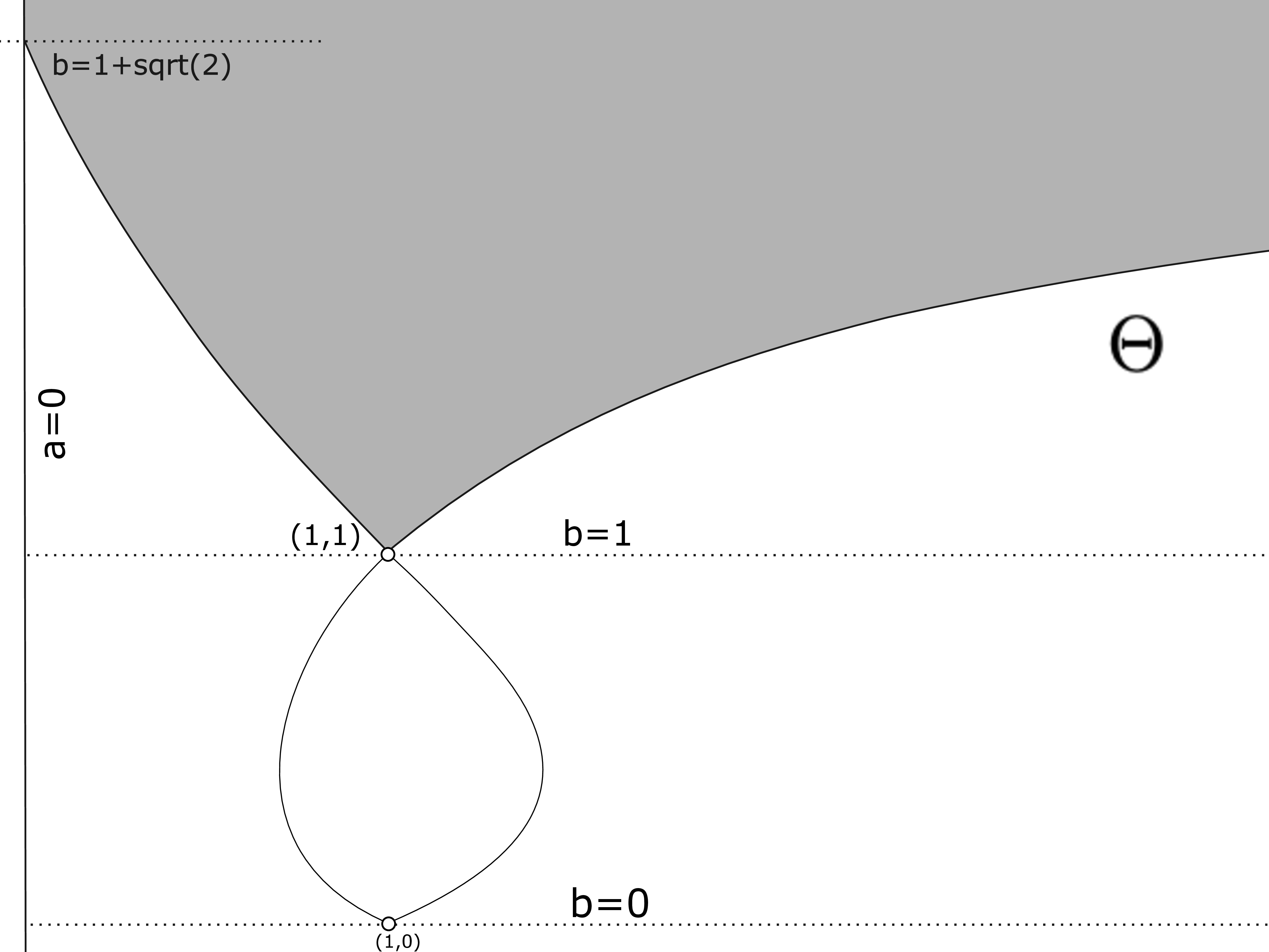}}
    \newline
    {\bf Figure 5.2\/}  The domain $\Theta$.
  \end{center}

  \subsection{Morphing the Groups}

  In this section we explain the representations
  constructed in [{\bf BLV\/}], but with our
  rational parametrization.  We take the same
  initial marked box $M_{c,d}$ shown in Figure 4.3.
  \newline
  \newline
  {\bf Modified Marked Box Operations:\/}
  As in [{\bf BLV\/}], we write $\lambda=(a,b)$.
    B-L-V define $3$ modified
      marked box operations. For each
        $\tau \in \{i,t,b\}$ they define
    \begin{equation}
      \tau^{\lambda}(M)=\sigma_{\lambda} \circ \tau(M).
    \end{equation}
    They show that these operations satisfy the
    same operations as the original ones and hence
    form a modular group of morphed marked box
    operations.     It turns out that this morphed
    marked box orbit still has a $\Z/3 * \Z/3$
    group of projective transformation symmetries.

   Given $(a,b) \in \Theta$ and $(c,d) \in (-1,1)^2$, the {\it morphed
   orbit\/}   is the orbit of $M_{c,d}$ under these morphed
   operations corresponding to $\lambda=(a,b)$.
   The boxes in this orbit are
       either disjoint or strictly nested.
    Using an argument akin to that in
    [{\bf S0\/}], B-L-V show that this property
    forces the corresponding representation of $\Z/3*\Z/3$ to
    be discrete and faithful.  Also, the strict nesting of the marked
    boxes forces the limit set to be a Cantor set.
    B-L-V  also show that their representations
    are {\it Anosov\/}.  See [{\bf BLV\/}] for definitions and the
    proof.  The central point is that most of the marked boxes
    in the orbit are small and thin.
        \newline
    \newline
    {\bf Order Two Symmetry:\/}
    The construction above gives a $4$-parameter family of
    representations
    of $\Z/3* \Z/3$.    B-L-V identify a
        certain function $h$ such that when
        $h(\lambda)=0$ there is a polarity $\sigma_2$ that
        conjugates the $\Z/3*\Z/3$ subgroup to
        itself, swapping the order $3$ element $\sigma_3$
        associated
        to $M$ and the order $3$ element associated to
        $i^{\lambda}(M)$.   B-L-V use an implicit function
        argument (which we re-do below) to show that
        for fixed $(c,d) \not = (0,0)$ the level curve
        $h(\lambda)=0$ is locally a smooth arc with endpoint
        $(a,b)=(1,1)$.  (This 
        is $(\epsilon,\delta)=(0,0)$ in their coordinates.)
       The group generated by $\sigma_2$ and $\sigma_3$ is
       the modular group representation associated to
       $(a,b,c,d)$.
       \newline
       \newline
       {\bf Remark:\/}
        Unlike in the Pappus case, this extra duality does not
        seem adapted to the morphed marked boxes.
        \newline
        \newline
        {\bf Explicit Formulas:\/}
        The two order $3$ matrices for the representation
        corresponding to $(a,b,c,d)$ are
        \begin{equation}
          r_1, \hskip 30 pt
          \Sigma_{a,b}^{-1} r_2 \Sigma_{a,b}.
        \end{equation}
        Here $r_1$ and $r_2$ are as in
        \S \ref{FORMULA}.  The matrix $r_1$ does not
        depend on the parameters $(a,b)$.
        Again, the product of these matrices has determinant $1$.
        My program checks that $r_1$ and $r_2$ have the
        same action as given in Equation \ref{action}, with
        $i^{(a,b)}$ replacing $i$, etc.

        \subsection{Extending the Representations}
    \label{DUALITY}

    Now we discuss the conditions on
    $(a,b,c,d)$ which guarantee that
    our representation of $\Z/3*\Z/3$ extends to
    a representation of $\Z/3*\Z/2$.
    \newline
    \newline
    {\bf The Duality Equation:\/}
As discussed in [{\bf BLV\/}] the necessary
and sufficient condition is that there is a polarity
conjugating $r_1$ to $r_2$.   We can compute the
condition in one of three ways:
\begin{enumerate}
\item $ {\rm det\/}(r_1r_2-I)=0$.  See
  [{\bf BLV\/}, Eq. 10.1].
\item Set $h(\epsilon,\delta)=0$ and change variables. See [{\bf BLV\/}] just
  after Eq. 10.1. 
\item ${\rm tr\/}(r_1r_2)-{\rm tr\/}(r_1^2r_2^2)=0$.  This is my formulation.
\end{enumerate}
Here ${\rm tr\/}$ is ``trace''.
Method 2 does not require us to compute the matrices
$r_1,r_2$ above.  Our code checks that the three methods give the same
equation.

For fixed $(c,d)$ we call the subset $\gamma_{c,d} \subset \Theta$
of parameters $(a,b)$ satisfying these conditions the
{\it duality curve\/}.
All the computations lead to the condition that
$\psi(a,b,c,d)=0$, where $\psi(a,b,c,d)$ is the following
expression.

{\small
\[
\begin{array}{rcl}
\psi(a,b,c,d) &=&
(a^2-1)(b^2+1)
\bigl(a^2b^2+a^2+ab^2-a+b^2+1\bigr)
\bigl(c^2+d^2-2c^2d^2\bigr)
\\[4pt] &&\quad
+\,a(b^2-1)
\bigl(a^2b^2+a^2+2ab^2-4ab-2a+b^2+1\bigr)
cd(c^2-d^2).
\end{array}
\]
\/}
(ChatGPT helped me find this way of writing it.)
Here $\psi$ is the numerator of a rational expression whose
denominator is $4a^2b^2(1-c^2)(1-d^2)$.
\newline
\newline
{\bf The Duality Curve:\/}
We will analyze $\gamma_{c,d}$ in the next chapter.
  One thing we note is that when
  $(c,d)=(0,0)$ we have $\psi=0$.  This
  case corresponds to the classic representations
  of the modular group which preserve a line
  in $\P$.   We treat this case specially.
  \newline
  \newline
  {\bf Symmetries:\/}
First,  $\psi \circ \theta_4=\psi$, where $\theta_4$ is as
in \S \ref{papspace}.  Second,
  $$({\rm tr\/}(r_1r_2)-{\rm tr\/}(r_1^2r_2^2))(1/a,b,d,c)=
 ( {\rm tr\/}(r^2_1r^2_2)-{\rm tr\/}(r_1r_2))(a,b,c,d).$$
  This implies what we call {\it Inverse Symmetry\/}:
  \begin{equation}
    {\rm sign\/}\ \psi(1/a,b,d,c)=-{\rm sign\/}\ \psi(a,b,c,d).
  \end{equation}
In particular,  $\psi(1/a,b,d,c)=0$ iff $\psi(a,b,c,d)=0$.
  \newline

  \noindent
  {\bf Local Calculation:\/}
  In [{\bf BLV\/}] the authors make a local calculation
  showing (in their coordinates) that when
  $(c,d) \not = (0,0)$ the set $\gamma_{c,d}$ is a
  smooth regular curve in a neighborhood of
  $(a,b)=(1,1)$.
  We make a similar calculation here.
  Define $\Phi: \R^2 \to \R^2$ by the map
\begin{equation}
  \Phi(a,b)=({\rm tr\/}(r_1r_2),{\rm tr\/}(r_1^2r_2^2)).
\end{equation}
We compute the Jacobian
\begin{equation}
\det(d\Phi) =
\frac{8(c^2+d^2-2c^2d^2)(c^2+d^2-2)}
{(1-c^2)^2(1-d^2)^2}.
\end{equation}
Since $(c,d) \not = (0,0)$, this determinant is nonzero.
Hence the set $\psi=0$ is a smooth curve in
a neighborhood of $(1,1)$.  The point
$(1,1)$ divides this curve into two arcs, and
the one corresponding to ${\rm trace\/}(r_1r_2)<-1$
is $\gamma_{c,d}$.  The other arc corresponds
to ${\rm tr\/}(r_1r_2) \in (0,1)$, and
here $r_1r_2$ is elliptic.
The corresponding
representation cannot be both discrete and
faithful. This shows that the
endpoint $(1,1)$ of $\gamma_{c,d}$, corresponding
to a Pappus group, lies in the boundary of
${\cal DFR\/}$ when ${\cal DFR\/}$
is considered as a subset of
$\cal R$.

This also shows, using the Implicit Function Theorem,
that the subset ${\cal P\/} \subset {\cal R\/}$,
consisting of the Pappus representations, is a smooth surface
away from the totally symmetric point.  The point here
is that the function ${\rm tr\/}(r_1r_2)$ is a smooth
function on $\cal R$ away from the symmetric point
and our matrix calculation shows that locally
this is a regular mapping on $\cal P$.

      \section{Proof of the Main Theorem}

  \subsection{Algebraic Tricks}

  Here we describe some algebraic tricks we use below.
  \newline
  \newline
  {\bf Resultants:\/}
The {\it resultant\/} of
    $P=a_2c^2+a_1c+a_0$ and
$Q=b_3c^3 + b_2 c^2 + b_1c + b_0$ is the number
    \begin{equation}
      {\rm res\/}(P,Q)=
      {\rm det\/} \left[\matrix{
          a_2 & a_1 & a_0 &0 &0 \cr
          0 & a_2 & a_1 & a_0 &0 \cr
          0& 0 & a_2 & a_1 & a_0 \cr
          b_3 & b_2 & b_1 &b_0 &0 \cr
          0 & b_3 & b_2 & b_1 &b_0}\right]
    \end{equation}
    The determinant vanishes if and only if $P$ and $Q$ have
    a common (complex) root.  The case for general pairs of
    polynomials works the same way;
    we just display the special case for
    typesetting purposes.
    See [{\bf Sil\/}, \S 2] for more.

    In the multivariable case, one can treat two polynomials
    $P(x_1,...,x_n)$ and $Q(x_1,...,x_n)$ as elements of the
    ring $R[x_n]$ where $R=\C[x_1,...,x_{n-1}]$.  The resultant
    ${\rm res\/}_{x_n}(P,Q)$
    computes the resultant in $R$ and thus gives
    a polynomial in $\C[x_1,...,x_{n-1}]$. The polynomials
    $P$ and $Q$ simultaneously vanish at
    $(x_1,...,x_n)$ only if
    ${\rm res\/}_{x_n}(P,Q)$ vanishes at
    $(x_1,...,x_{n-1})$.
    \newline
    \newline
    {\bf Taylor series:\/}
Given a polynomial
$H(...,b,...)$ we define
\begin{equation}
  H^{(k)} = \frac{\partial^k H}{\partial b^k}\bigg|_{b=1}.
\end{equation}
To prove that $H(...b...) \geq \lambda$ it suffices to show
  \begin{itemize}
  \item  $H^{(0)} \geq \lambda.$
  \item $H^{(k)}(...)>0$ for all
    relevant variables and all $k \leq m$.
   \item $H^{(k)}=0$ for all $k>m$.
   \end{itemize}
We call this the
{\it Taylor series method\/} for obvious reasons.
When we use this below, the variables and domain will be clear.
\newline
\newline
{\bf Special Polynomials:\/}
A certain family of polynomials arises repeatedly in our analysis
below.
We treat this family here, in isolation.

\begin{lemma}
         \label{specialp}
         If $|\lambda| \leq 1$ then
         $f(c,d)=c^2+d^2-2c^2d^2 + \lambda( c^3d-cd^3)>0$
        on the domain $(-1,1)^2-\{0,0\}$.
      \end{lemma}

    \startproof
    Since $f(c,d)=f(-d,c)=f(-c,-d)=f(d,-c)$, 
    it suffices to prove this for  $ c,d \geq 0$.
    The result is obvious if $cd=0$, so we take $c,d \in (0,1)$.

    Define $A(u,v)=(u^2+v^2)-f(u,v)$.  Here $A$ is homogeneous of
    degree $4$.  Also, $f(u,v)=u^2+v^2 - A(u,v)$.
    If $u,v \in [0,1]$ and
   $\max(u,v)=1$ we have $f(u,v) \geq 0$,
    because
    $$f(u,1)=(1-u^2)(1- \lambda u) \geq 0 \hskip 30 pt
    f(1,v)=(1-v^2)(1+\lambda v) \geq 0.
    $$
    Hence, for these choices of $u$ and $v$, we have
    $A(u,v) \leq u^2+v^2$.
    Since $u^2+v^2>0$ we have
    $r^2A(u,v)<u^2+v^2$ when $r \in (0,1)$.
    When $c,d \in (0,1)$ we can write
    $(c,d)=(ru,rv)$ where $r \in (0,1)$ and $u,v \in [0,1]$ and $\max(u,v)=1$.
    We compute
    $$r^{-2} f(c,d)=u^2+v^2 - r^2A(u,v)>0.\spadesuit$$

    \subsection{The Walls of the Good Parameter Set}

    Let $S_b$ denote the
    set of all $a$ such that $(a,b) \in \Theta$.
    This is either a horizontal segment or ray,
    depending on $b$. Lemma \ref{GOOD} describes $S_b$.
        
  \begin{lemma}
    \label{WALL}
    For fixed $(c,d) \not = (0,0)$, and
    any $b>1$ we have $\psi<0$ at the
    left endpoint of $S_b$ and $\psi>0$ at
    the right endpoint of $S_b$ when it is finite,
    and otherwise for all sufficiently large values of $a$.
  \end{lemma}

  \startproof
  We first take $b \geq 1+\sqrt 2$, so
  that $S_b=(0,\infty)$.     We have
    \begin{equation}
      \label{cushion1}
      \psi(0,b,c,d)=(1+b^2)^2(2c^2d^2 - c^2-d^2)<0.
    \end{equation}
    By Inversive Symmetry,
    $\psi(a,b,c,d)>0$ when $a$ is sufficiently large.

    Now we take $b \in (1,1+\sqrt 2)$.
    The left endpoint of $S_b$ is
    \begin{equation}
      \label{constraint}
  a=\frac{1+2b-b^2}{1+b^2}.
\end{equation}
Making this substitution, we get
  $$
  \psi=-\bigg(\frac{4b(b^2-1)}{(1+b^2)^2}\bigg) \times \mu,
  $$
Now we compute
\[
\begin{array}{rcl}
\mu^{(0)}&=&8c^2+8d^2-16c^2d^2\\[3pt]
\mu^{(1)}&=&12c^2+12d^2-24c^2d^2-4c^3d+4cd^3\\[3pt]
\mu^{(2)}&=&12c^2+12d^2-24c^2d^2-4c^3d+4cd^3\\[3pt]
\mu^{(3)}&=&12c^2+12d^2-24c^2d^2+12c^3d-12cd^3\\[3pt]
\mu^{(4)}&=&24c^2+24d^2-48c^2d^2+24c^3d-24cd^3.
\end{array}
\]

All higher derivatives vanish.
We have $\mu>0$ by Lemma \ref{specialp} and the Taylor series method.
Hence $\psi<0$ at the left endpoint of $S_b$.
By Inversive Symmetry, we have
$\psi>0$ on the right endpoint.
\endproof

  \subsection{Counting the Zeros}

  \begin{lemma}
    \label{ZERO}
    If $(c,d) \not = (0,0)$ then $\psi=0$ exactly once in $S_b$.
  \end{lemma}

  \startproof
  Let $P_b(a)=\psi(a,b,c,d)$.
  For $b$ near $1$, our local calculation shows that
  $P_b$ has exactly one root in $\Theta$.
  If this situation ever changes as $b$ increases,
  there will be some value of $b$ such that
  $P_b(a)$ has a double root.  We rule this out by showing that
  $P_b$ and $dP_b/da$ never vanish simultaneously.  We compute
$$
{\rm res\/}_c(P_b,dP_b/da)=
    4(b^4-1)^3(d^2-1)^2 d^9 \times r(a,b)^3,
    $$
\begin{equation}
  \label{res}
\begin{array}{rcl}
r(a,b) &=& (a^6+1)(b^2+1)^2
 +(a^5+a)(4b^4-8b^3-8b-4) \\[3pt]
&& {}+(a^4+a^2)(5b^4-4b^3+2b^2+4b+5)
 +a^3(4b^4-4).
\end{array}
\end{equation}

  We just have to see that $r(a,b)$ does not
  vanish on $\Theta$.
  We first take $b \in (1,1+\sqrt 2)$ and restrict $r(a,b)$ to the
  two constraint curves $\alpha$ and $\beta$ given by
  Lemma \ref{GOOD}.  Equation \ref{constraint} gives the equation
  for $\alpha$ and we invert this equation to get the equation for
  $\beta$.
  When we do the restricting,
    we get
  $$r|_{\alpha}=\frac{8b(b^2-1)P(b)}{(1+b^2)^4},
  \hskip 30 pt
  r|_{\beta}=\frac{8b(b^2-1)P(b)}{(b^2-2b-1)^6},$$
\[
\begin{array}{rcl}
P(b) &=& 32+48(b-1)+24(b-1)^2+56(b-1)^3+92(b-1)^4 \\[3pt]
&&{}+52(b-1)^5+10(b-1)^6+2(b-1)^7+(b-1)^8.
\end{array}
\]
     Hence $r$ is positive on $\alpha \cup \beta$.
     
Now, every other point of $\Theta$ can be reached from a point
on $\alpha \cup \beta$ by an upward  vertical path.
So, to finish the proof,
it suffices to prove that
$\partial r/\partial b>0$ on $(0,\infty) \times (1,\infty)$.
The Taylor series method shows this:
\[
\begin{array}{rcl}
r^{(1)}(a) &=& (8-16a+16a^2)+16a^3+8a^4(2-2a+a^2)>0, \\[4pt]
r^{(2)}(a) &=& 16(1+a^6)+40(a^2+a^4)+48a^3>0, \\[4pt]
r^{(3)}(a) &=& 24(1+a^6)+48(a+a^5)+96(a^2+a^4)+96a^3>0, \\[4pt]
r^{(4)}(a) &=& 24(1+a^6)+96(a+a^5)+120(a^2+a^4)+96a^3>0.
\end{array}
\]
All higher derivatives vanish.
\endproof

\subsection{Geometry of the Duality Curves}

In this section we get some bounds on the duality curve.

\begin{lemma}
  \label{BOUND1}
  $\gamma_{c,d} \subset [1,2] \times [1,\infty)$
  if $0 \leq c \leq d$.
\end{lemma}

\startproof
For fixed $b,c,d$, the polynomial $\psi$ vanishes once.
We compute
$$\psi(1,b,c,d)=4b(b-1)^2(b+1)cd(c^2-d^2) \leq 0,$$
\[
\begin{array}{rcl}
\psi(2,b,c,d) &=&
 +16b^2(b-1)cd(d^2-c^2) \\[3pt]
&&{}+ \big(9c^2 + 9d^2 - 18c^2d^2 + (2cd^3 - 2c^3d)\big) \\[3pt]
&&{}+ \bigl(30c^2+30d^2-60c^2d^2+(16cd^3-16c^3d)\bigr)b^2 \\[3pt]
&&{}+ \bigl(21c^2+21d^2-42c^2d^2-(18cd^3+18c^3d)\bigr)b^4 \geq 0 .
\end{array}
\]
We used
Lemma \ref{specialp} three times here.
By the Intermediate Value Theorem,
$\psi$ vanishes somewhere on
$[1,2] \times \{b\}$.
\endproof

\noindent
{\bf Remark:\/} It follows from symmetry that
$\gamma_{c,d} \subset [1/2,2] \times [1,\infty)$ in
all cases.

\begin{lemma}
  As $(c,d) \to (0,0)$ the curve
  $\gamma_{c,d}$ converges to
  $\{1\} \times [1,\infty)$ uniformly
  on compact subsets.
\end{lemma}

\startproof
Assume that $\|(c,d)\|=\epsilon$ and we take $b \in [1,B]$
for any bound $B$.
We will show that the restriction
of $\gamma_{c,d}$ to the set $(0,\infty) \times [1,B]$
converges uniformly to $\gamma_{0,0}$ as
$\epsilon \to 0$.
We compute
\begin{equation}
  \label{small}
  \psi(1+t,b,c,d)=4b(b+1)(b-1)^2cd(c^2-d^2)+
  2(c^2+d^2)t + E_1t + E_2t^2
\end{equation}
where $E_1$ is a polynomial in which every term
has total degree at least $4$ in $c,d$ and
$E_2$ is a polynomial in which every term
has total degree at least $2$ in $c,d$.
But then, once $\epsilon$ is small enough,
$\psi(1,b,c,d)=O(\epsilon^4)$ and
$v(t)=\psi(1+t,b,c,d)$ varies by $O(\epsilon^3)$ on
each interval $[-\epsilon,0]$ and
$[0,\epsilon]$ and hence vanishes in
one of these intervals.
\endproof

Now we deal with the case $(c,d)=(0,0)$.
In this case $\psi=0$ identically.  However,
in this case we have a redundant description of
our representations.  To remove the redundancy,
we set $a=1$.  The trace of $r_1r_2$ is
$-(3b^2-1)^2/4b^2$.
As $b$ varies in $(1,\infty)$ we get every trace
in $(-\infty,-1)$.  Thus, if we {\it define\/}
$\gamma_{0,0}$ to be the vertical ray
$\{1\} \times [1,\infty)$, we pick up all the conjugacy
classes of the representations.  This definition
makes $\gamma_{c,d}$ vary continuously with $(c,d)$.

\subsection{The Proof modulo Properness}

Now we know that $\gamma_{c,d}$ is a
curve that starts at $(1,1)$, intersects every horizontal
segment $S_b$ exactly once, and stays
$[1/2,2] \times [1,\infty)$.

Let $\cal C$  be the space of Pappus representations
as in \S \ref{papspace}.
 We introduce a new space $\cal H$.  This space is
 a fiber bundle over the space ${\cal C\/}$.
The fiber  over $[(c,d)]$ is the set of representations corresponding to
$\gamma_{c,d}$.   This makes saense because
  the representations on 
  $\gamma_{-d,c}$ are conjugate to those on
  $\gamma_{c,d}$.
  The space $\cal H$ is homeomorphic to the
  upper half-space and the representations vary
  continuously.  That is, we have a
  continuous map $f: {\cal H\/} \to {\cal R\/}$,
  the representation space.  On $\cal C$ the map
  $f$ is the map considered in \S \ref{papspace},
  and we have ${\cal P\/}=f({\cal C\/})$.
  
  We let $\cal B$ be the component of ${\cal R\/}-{\cal P\/}$
  which does not contain the origin. The other side of
  $\cal B$, near $\cal P$, consists of representations that are not
  both discrete and faithful.  Hence $f({\cal H\/}) \subset \cal
  B$.   Below we will prove that $f$ is a proper map.
  We let $\widehat {\cal H\/}$ be the
  $1$-point compactification of $\cal H$.
  Likewise, we let $\widehat {\cal B\/}$ denote
  the $1$-point compactification of $\cal B$.
  Both spaces are homeomorphic to closed $3$-balls.
  Since $f$ is proper, $f$
  extends to a map from
  $\widehat {\cal H\/}$ to $\widehat {\cal B\/}$ which
  is a homeomorphicm from the sphere
  ${\cal C\/} \cup \{\infty\}$ to the sphere
  ${\cal P\/} \cup \infty$.
  But then $f$ is surjective.  That is,
  $f({\cal H\/})={\cal B\/}$.
  Hence $\cal B$ is precisely the
  component of $\cal DFR$ that contains
  $\cal P$.
  
  This completes the proof of Theorem \ref{one}
  modulo the statement that $f$ is a proper map.

\subsection{Proof of Properness}
  
Now we prove that $f$ is proper.
Rather than work directly with $\cal H$ we
pass to the $4$-fold branched cover,
which is  the subset of points $(a,b,c,d)$ in
$\Theta \times (-1,1)^2$ such that $(a,b) \in \gamma_{c,d}$.
We suppose that we have a sequence
$\{(a_n,b_n,c_n,d_n)\}$ of parameters
that exits every compact subset of our
domain.  We want to see that in all cases
the corresponding representation exits
every compact subset of $\cal R$.
By symmetry, it suffices to consider the
case when $c_n,d_n\geq 0$.

To avoid repetitive calculations we assume
that $c_n \leq d_n$.  This gives $a_n \in [1,2]$.
The other case is
similar.  In Cases 2 and 3 below, for the other case, we
would use the
$r_1^2r_2$ in place of the element $r_1r_2^2$ below
and we would have $a_n \in [1/2,1]$ rather
than $a_n \in [1,2]$.

After passing to a subsequence we arrive
at $3$ cases which cover everything:
\begin{enumerate}
\item
  $a_n \to a \in [1,2]$ and
  $ b_n \to \infty$ and
  $c_n \to c \in [0,1]$ and
  $d_n \to d \in [0,1]$.

\item $a_n \to a \in [1,2]$
  and $b_n \to b \in [1,\infty)$
  and $c_n \to c \in [0,1)$
  and $d_n \to 1$.

\item $a_n \to a \in [1,2]$
  and $b_n \to b \in [1,\infty)$
  and $c_n \to 1$ and $d_n \to 1$.
\end{enumerate}
Again, in all these cases we have $0 \leq c_n \leq d_n<1$.
It suffices in all cases to show that the trace of some word,
when it is normalized to have unit determinant,
tends to $\infty$ with $n$.  Equivalently -- and without
normalizing the determinant -- it suffices
to show that the conjugacy invariant in Equation \ref{limit}
tends to $\infty$ for a suitable word.
       \newline
  \newline
  {\bf Case 1:\/}
  We use $m=r_1r_2$.   We compute
  \begin{equation}
    {\rm trace\/}(m_n)=
    \frac{-U_nV_nb_n^4+W_n}{4a_n^2b_n^2(1-c_n^2)(1-d_n^2)} \sim
    \frac{-U_nV_nb_n^2}{4a_n^2(1-c_n^2)(1-d_n^2)},
  \end{equation}
\[
U_n = a_n^2(1-c_n^2)+ (1+a_n)(1-c_nd_n)>a_n^2(1-c_n^2).
\]
\[
  V_n=(1+a_n)(1+c_nd_n) +a_n^2(1-d_n^2)>2.
\]
$W_n$ is a polynomial in which $b_n$ appears with maximum degree $3$.
From all this, we see
easily that ${\rm trace\/}(m_n) \to -\infty$ as
$n \to \infty$.
\newline
\newline
{\bf Case 2:\/}
We use $m=r_1r_2^2$ and $\tau(m)$ as
in Equation \ref{TAU}.
We compute

\begin{equation}
  \label{limit1}
  \tau(m_n):=\frac{Y(a_n,b_n,c_n,d_n)^3}{64 a_n^6b_n^6(1-c_n^2)^4(1-d_n^2)^2}.
\end{equation}
We compute that $Y(a,b,c,1)=(1-c^2)Z$, where
{\footnotesize
\[
\begin{array}{rcl}
Z^{(0)} &=&
4(1+2a^2+a^4)+4c(a-1)+4a^3c(a-1) \geq 16
\\[4pt]

Z^{(1)} &=&
8(1-c)+8a+24a^2+8a^3+8a^4
   +4ac+4a^3c(2a-1)>0
\\[4pt]

Z^{(2)} &=&
16+24a+8a^2(8-c^2)+24a^3+16a^4
   +16c(a^4-1)>0
\\[4pt]

Z^{(3)} &=&
24+48a+24a^2(4-c^2)+48a^3+24a^4
   +24c(a^4-1)+12ac(a^2-1)>0
\\[4pt]

Z^{(4)} &=&
24+48a+24a^2(3-c^2)+48a^3+24a^4
   +24c(a^4-1)+24ac(a^2-1)>0
\end{array}
\]
\/}
All higher derivatives vanish.
Hence $Z \geq 16$, by the Taylor series method.
Hence $Y^3(a,b,c,1) \geq 4096(1-c^2)^3>0$.  This bound
combines with Equation \ref{limit1} to show that $\tau(m_n) \to
\infty$.
\newline
\newline
{\bf Case 3:\/}
We keep the same notation from Case 2.
We have $Y(a,b,1,1)=0$. We expand
$Y$ in a Taylor series about $(c,d)=(1,1)$ and find that
\begin{equation}
  Y(a,b,c,d)=U(1-c)+ V (1-d) + {\rm higher\ order\ terms\/}.
\end{equation}
Here $$U=-\frac{\partial Y}{\partial c}\bigg|_{(c,d)=(1,1)}, \hskip 40 pt
V=-\frac{\partial Y}{\partial d}\bigg|_{(c,d)=(1,1)}.$$
Let $W_{\pm}=U \pm V$.  The Taylor series method
shows that $W_{\pm} \geq 32$ on $[1,2] \times [1,\infty)$:
{\small
\[
\begin{array}{rcl@{\hskip 30pt}rcl}
W_+^{(0)}&=&32a^2-16a^3+16a^4 \geq 32
&
W_-^{(0)}&=&16a+16a^4 \geq 32 \\[3pt]
W_+^{(1)}&=&16a+64a^2+32a^4
&
W_-^{(1)}&=&32a+32a^2+16a^3+32a^4\\[3pt]
W_+^{(2)}&=&48a+128a^2+48a^3+64a^4
&
W_-^{(2)}&=&48a+96a^2+48a^3+64a^4\\[3pt]
W_+^{(3)}&=&96a+192a^2+144a^3+96a^4
&
W_-^{(3)}&=&48a+96a^2+96a^3+96a^4\\[3pt]
W_+^{(4)}&=&96a+192a^2+192a^3+96a^4
&
W_-^{(4)}&=&96a^3+96a^4.
\end{array}
\]
\/}
Since $U \pm V \geq 32$ we also have
$U \geq 32$.
We set $\epsilon_n=(1-c_n)$.  Since $c_n \leq d_n$ we have
$\epsilon_n\geq 1-d_n>0$.  This gives
$$Y(a_n,b_n,c_n,d_n) \sim U \epsilon_n + V (1-d_n)
\geq \min(U,U+V) \epsilon_n \geq 32\epsilon_n.$$
Hence $Y(a_n,b_n,c_n,d_n)>\epsilon_n$
once $n$ is sufficiently large.
Hence, the numerator in Equation \ref{limit1} is at least
$\epsilon_n^3$ for $n$ large.
But the denominator in Equation \ref{limit1} is at most
$2^{12}a_n^6b_n^6 \epsilon_n^6 \sim 2^{12}a^6b^6 \epsilon_n^6$.
Hence the whole
expression tends to $\infty$ as $n \to \infty$.
\newline

This completes the proof that the map
$f: {\cal H\/} \to {\cal R\/}$ is proper.

\newpage
      \section{References}

\noindent
[{\bf Bar\/}]
T. Barbot,
{\it Three dimensional Anosov Flag Manifolds\/},  Geometry $\&$ Topology (2010)
\vskip 8 pt
\noindent
[{\bf BLV\/}], T. Barbot, G. Lee, V. P. Valerio, {\it Pappus's Theorem,
  Schwartz Representations, and Anosov Representations\/},
Ann. Inst. Fourier (Grenoble) {\bf 68\/} (2018) no. 6
\vskip 8 pt
\noindent
[{\bf BCLS\/}] M. Bridgeman, D. Canary, F. Labourie, A. Samburino,
{\it The pressure metric for Anosov representations\/}, GAFA  {\bf 25\/} (2015)
\vskip 8 pt
\noindent
[{\bf FL\/}] C. Florentino, S. Lawton, {\it The topology of moduli
  spaces of free groups\/}, Math. Annalen {\bf 345\/}, Issue 2 (2009)
\vskip 8 pt
\noindent
[{\bf GW\/}] O. Guichard, A. Wienhard, {\it Anosov Representations: Domains of
  Discontinuity and applications\/}, Invent Math {\bf 190\/} (2012)
\vskip 8 pt
\noindent
[{\bf Hit\/}] N. Hitchin, {\it Lie Groups and Teichmuller Space\/}, Topology {\bf 31} (1992)
\vskip 8 pt
\noindent
[{\bf KL\/}] M. Kapovich, B. Leeb, {\it Relativizing characterizations of
  Anosov subgroups, I (appendix by Gregory A. Soifer).\/} Groups Geom. Dyn.  {\bf 17\/}  (2023)
  \vskip 8 pt
  \noindent
[{\bf Lab\/}] 
F. Labourie, {\it Anosov Flows, Surface Groups and Curves in Projective Spaces\/},
P.A.M.Q {\bf 3\/} (2007)
\vskip 8 pt
\noindent
[{\bf L\/}] S. Lawton, {\it Generators, relations, and symmeries in
  pairs of $3 \times 3$ unimodular matrices\/}, J. Algebra, 313(2)
(2007)
\vskip 8 pt
\noindent
[{\bf S0\/}] R. E. Schwartz, {\it Pappus's Theorem and the Modular Group\/},
Publ. IHES (1993)
\vskip 8 pt
\noindent
[{\bf S1\/}] R. E. Schwartz, {\it Le Retour de Pappus\/},
KIAS-Springer Lecture Notes (2025) to appear.
See also arXiv 2412.02417
\vskip 8 pt
\noindent
[{\bf Sil\/}],
  J. Silverman, {\it The arithmetic of dynamical systems\/}, Graduate
  Texts in Mathemtics {\bf 241\/} (2007) Springer
\vskip 8 pt
\noindent
[{\bf V\/}], V. P. Valerio, {\it Teorema de Pappus, Representa\c{c}oes de Schwartz e
  Representa\c{c}oes Anosov\/},
Ph. D. Thesis, Federal University of Minas Gerais (2016)
\vskip 8 pt
\noindent
[{\bf W\/}] S. Wolfram et. al., {\it Mathematica\/}, Version 11
Wolfram Research Inc. (2024)

\end{document}